\documentclass[11pt,amssym,twoside]{article}
\usepackage{amssymb}
\newtheorem{thm}{Theorem}[section]

\newtheorem{conj}[thm]{Conjecture}

\newcommand{\To}{\longrightarrow}

\newcommand{\Z}{\mathbb Z}

\newcommand{\Q}{\mathbb Q}
\newcommand{\LL}{\mathcal {L}}
\newcommand{\N}{\mathcal {N}}
\newcommand{\gal}{\mathcal {G}al}
\newcommand{\out}{\mathcal {O}ut}
\begin{document}

\title{Deformation of Outer Representations of Galois Group II}
\author{Arash Rastegar}


\maketitle
\begin{abstract}
This paper is devoted to deformation theory of "anabelian"
representations of the absolute Galois group landing in outer
automorphism group of the algebraic fundamental group of a
hyperbolic smooth curve defined over a number-field. In the first
part of this paper, we obtained universal deformations for
Lie-algebra versions of the above representation using the
Schlessinger criteria for functors on Artin local rings. In the
second part, we use a version of Schlessinger criteria for
functors on the Artinian category of nilpotent Lie algebras which
is formulated by Pridham, and explore arithmetic applications.
\end{abstract}
\section*{Introduction}

Based on Grothendieck's anabelian philosophy, the Galois-module
structure of outer automorphism group of the geometric
fundamental group of a smooth hyperbolic curve should contain all
the arithmetic information coming from $X$. Particularly, the so
called "Hom" conbjecture, which is one of Grothendieck's anabelian
conjectures, proved by Mochizuki [Moc], states that, for smooth
$X$ and $X'$ hyperbolic curves defined over a number field $K$,
there is a natural one-to-one correspondence
$$
Isom_{K}(X,X')\To Out_{Gal(\bar {K}/K)}(Out(\pi_1(\bar { X}))
,Out(\pi_1(\bar {X}'))).
$$
Here $\pi_1(\bar {X})$ denotes the geometric fundamental group of
$X$ and $Out(\pi_1(\bar {X}))$ denotes its outer automorphism
group, i.e. the quotient of the automorphism group $Aut(\pi_1(\bar
{X}))$ by inner automorphisms. By $Out_{Gal(\bar {K}/K)}$ we mean
the set of Galois equivariant isomorphisms between the two
profinite groups divided by the inner action of the second
component. The Galois structure of $Out(\pi_1(\bar {X}))$ is given
by a continuous group homomorphism
$$
\rho_X:Gal(\bar {K}/K) \To Out(\pi_1(\bar {X}))
$$
associated to the following short exact sequence
$$
0 \To \pi_1(\bar {X}) \To \pi_1(X) \To Gal(\bar {K}/K) \To 0.
$$
By a result of Grothendieck, the outer representation breaks to an
outer representation of $Gal(K_S/K)$ for a finite set $S$ of
places of $K$, where $K_S$ is the maximal extension of $K$
unramified outside $S$.

One can induce a filtration on the Galois group using the induced
pro-$l$ representation
$$
\rho^l_X:Gal(K_S/K) \To Out(\pi_1(\bar X)^{(l)})
$$
and the weight filtration on the outer automorphism group. These
filtrations induce graded nilpotent $\mathbb {Z}_l$-Lie algebras
on both sides and a representation of graded Lie algebras
$$
\gal (K_S/K) \To \out (\pi_1(\bar X)^{(l)})
$$
which is the object we plan to deform.

In the first part of this paper, we used the classical
Schlessinger criteria for deformations of functors on Artin local
rings for deformation of the above representation [Sch]. Here, we
use a graded version of Pridham's adaptation of Schlessinger
criteria for functors on finite dimensional nilpotent Lie
algebras [Pri]. The ultimate goal is to have the whole theory of
deformations of Galois representations and its relation to
modular forms, translated to the language of nilpotent Lie
algebras and their representations and use its computational and
conceptual advantages. It is also natural to formulate a Lie
algebra version of Grothendieck's anabelian conjectures.

\section{Lie algebras associated to profinite groups}

Exponential map on the tangent space of an algebraic group
defined over a field $k$ of characteristic zero, gives an
equivalence of categories between nilpotent Lie algebras of
finite dimension over $k$ and unipotent algebraic groups over $k$.
This way, one can associate a Lie algebra to the algebraic
unipotent completion $\Gamma^{alg}(\Q)$ of any profinite group
$\Gamma$.

On the other hand, Mal\v{c}ev defines an equivalence of
categories between nilpotent Lie algebras over $\Q$ and uniquely
dividible nilpotent groups. Inclusion of such groups in nilpotent
groups has a right adjoint $\Gamma\to \Gamma_{\Q}$. For a nipotent
group $\Gamma$, torsion elements form a subgroup $T$ and
$\Gamma_{\Q}= \cup (\Gamma/T)^{1/n}$. In fact we have
$\Gamma_{\Q}=\Gamma^{alg}(\Q)$.

Any nilpotent finite group is a product of its sylow subgroups.
Therefore, the profinite completion $\Gamma^{\wedge}$ factors to
pro-$l$ completions $\Gamma_l^{\wedge}$, each one a compact open
subgroup of the corresponding $\Gamma^{alg}(\Q_l)$ and we have
the following isomorphisms of $l$-adic Lie groups
$$
\mathrm { Lie}(\Gamma_l^{\wedge})=\mathrm {Lie}
(\Gamma^{alg}(\Q_l))=\mathrm {Lie}(\Gamma^{alg}(\Q))\otimes \Q_l.
$$
In fact, the adelic Lie group associated to $\Gamma$ can be
defined as $\mathrm {Lie}(\Gamma^{alg}(\Q))\otimes \mathbb {A}^f$
which is the same as $\prod \mathrm {Lie}(\Gamma_l^{\wedge})$.

Suppose we are given a nilpotent representation of $\Gamma$ on a
finite dimensional vector space $V$ over $k$, which means that
for a filtration $F$ on $V$ respecting the action, the induced
action on $Gr_F(V)$ is trivial. The subgroup
$$
\{\sigma\in GL(V)|\sigma F=F , Gr_F(\sigma)=1\}
$$
is a uniquely divisible group and
one obtains a morphism
$$
\mathrm {Lie}(\Gamma_{\Q})\to \{\sigma\in gl(V)| \sigma F=F ,
Gr_F(\sigma)=0\}
$$
which induces an equivalence of categories between nilpotent
representations of $\Gamma$ and nilpotent representations of the
Lie algebra $\mathrm {Lie}(\Gamma_{\Q})$ over the field $k$. The
above equivalence of categories extends to an equivalence between
linear representations of $\Gamma$ and representations of its
algebraic envelope [Del].

The notion of weighted completion of a group developed by Hain
and Matsumoto generalizes the concept of algebraic unipotent
completion [Hai-Mat]. Suppose that $R$ is an algebraic $k$-group
and $w:\mathbb G_m\to R$ is a central cocharacter. Let $G$ be an
extension of $R$ by a unipotent group $U$ in the category of
algebraic $k$-groups
$$
0\To U\To G\To R\To 0.
$$

The first homology of $U$ is an $R$-module, and therefore a
$\mathbb G_m$-module via $w$, which naturally decomposes to
direct sum of irreducible representations each isomorphic to a
power of the standard character. We say that our extension is
negatively weighted if only negative powers of the standard
character appear in $H_1(U)$. The weighted completion of $\Gamma$
with respect to the representation $\rho$ with Zariski dense image
$\rho:\Gamma\to R(\Q_l)$ is the universal $\Q_l$-proalgebraic
group $\mathcal G$ which is a negatively weighted extension of $R$
by a prounipotent group $\mathcal U$ and a continuous lift of
$\rho$ to $\mathcal G(\Q_l)$ [Hai-Mat]. The Lie algebra of
$\mathcal G(\Q_l)$ is a more sophisticated version of
$\mathrm{Lie}(\Gamma_{\Q})\otimes \Q_l$.

\section{Functors on nilpotent graded Lie algebras}

In this section, we review Pridham's nilpotent Lie algebra version
of Schlessinger criteria [Pri]. The only change we impose is to
consider finitely generated graded nilpotent Lie algebras with
finite dimensional graded pieces, instead of finite dimensional
nilpotent Lie algebras.

Fix a field $k$ and let $\N_k$ denote the category of finitely
generated NGLAs (nilpotent graded Lie algebras) with finite
dimensional graded pieces, and $\widehat {\N_k}$ denote the
category of pro-NGLAs with finite dimensional graded pieces which
are finite dimensional in the sense that $\dim L/[L,L]<\infty$.
Given $\LL \in \widehat {\N_k}$ define $\N_{\LL,k}$ to be the
category of pairs $\{ N\in \N_k, \phi :\LL\to N\}$ and $\widehat
{\N_{\LL,k}}$ to be the category of pairs $\{ N\in \widehat
{\N_k}, \phi :\LL\to N\}$.

All functors on $\N_{\LL,k}$ should take the $0$ object to a one
point set. for a functor $F:\N_{\LL,k}\to \textrm {Set}$, define
$\hat F:\widehat {\N_{\LL,k}}\to \textrm {Set}$ by
$$
\hat F(L)=\lim_{\leftarrow} F(L/\Gamma_n(L)),
$$
where $\Gamma_n(L)$ is the $n$-th term in the central series of
$L$. Then, for $h_L:\N_{\LL,k}\to \textrm {Set}$ defined by $N\to
Hom(L,N)$ we have an isomorphism
$$
\hat F(L) \longrightarrow Hom(h_L,F)
$$
which can be used to define the notion of a pro-representable
functor.

A morphism $p\in N\to M$ in $\N_{\LL,k}$ is called a small
section if it is surjective with a principal ideal kernel $(t)$
such that $[N,(t)]=(0)$.

Given $F:\N_{\LL,k}\to \textrm {Set}$, and morphisms $N'\to N$
and $N''\to N$ in $\N_{\LL,k}$, consider the map
$$
F(N'\times_N N'')\To F(N')\times_{F(N)} F(N'').
$$
Then, by the Lie algebra analogue of the Schlessinger theorem $F$
has a hull if and only if it satisfies the following properties

(H1) The above map is surjective whenever $N''\to N$ is a small
section.

(H2) The above map is bijective when $N=0$ and $N''=L(\epsilon )$.

(H3) $\dim_k (t_F)< \infty$. \\ $F$ is pro-representable if and
only if it satisfies the following additional property

(H4) The above map is an isomorphism for any small extension
$N''\to N$. \\ Note that, in case we are considering graded
deformations of graded Lie algebras, only the zero grade piece of
the cohomology representing the tangent space shall be checked to
be finite dimensional.

\section{Several deformation problems}

Let $X$ denote a hyperbolic smooth algebraic curve defined over a
number field $K$. Let $S$ denote the set of bad reductions of $X$
together with places above $l$. We shall construct Lie algebra
versions of the pro-$l$ outer representation of the Galois group
$$
\rho^l_X:Gal(K_S/K) \To Out(\pi_1(\bar X)^{(l)}).
$$
Let $I_l$ denote the decreasing filtration on
$Out(\pi_1(X)^{(l)})$ induced by the weight filtration of
$\pi_1(\bar X)^{(l)}$. By abuse of notation, we also denote the
filtration on $Gal(K_S/K)$ by $I_l$. We get an injection of the
associated graded $\Z_l$-Lie algebras on both sides
$$
\gal (K_S/K) \To \out (\pi_1(\bar X)^{(l)}).
$$
One can also start with the $l$-adic unipotent completion of the
fundamental group and the outer representation of Galois group on
this group.
$$
\rho^{un,l}_X:Gal(K_S/K) \To Out(\pi_1(\bar X)^{un}_{/\Q_l}).
$$
and the associated Galois Lie algebra would be the same as those
associated to $I_l$ ([Hai-Mat] 8.2). Let $U_S$ denote the
prounipotent radical of the zariski closure of the image of
$\rho^{un,l}_X$. The image of $Gal(K_S/K)$ in $Out(\pi_1(\bar
X)^{un}_{/\Q_l})$ is a negatively weighted extension of $\mathbb
G_m$ by $U_S$ with respect to the central cocharacter $w:x\mapsto
x^{-2}$. The weight filtration induces a graded Lie algebra
$\mathcal U_S$ which is isomorphic to $\gal (K_S/K)\otimes \Q_l$
([Hai-Mat] 8.4).

There are several deformation problems in this setting which are
interesting. For example, the action of Galois group on unipotent
completion of the fundamental group induces an action of the
Galois group on the corresponding nilpotent $\Q_l$-Lie algebra
$$
\rho^{ni,l}_X:Gal(K_S/K) \To Aut(\mathcal U_S)
$$
which could be deformed using the Pridham's version of
Schlessinger criteria. We will explain in the following section,
why this representation is completely determind by the
abelianized representation of the Galois group. We will use
results of Koneko mentioned in the first part of the paper.
Therefore, deformation theory of this object is the same as the
abelianized deformation theory. Although, in this formulation we
get universal deformation nilpotent Lie algebras instead of
univeral deformation rings.

A similar thing to do would be deforming the following
representation
$$
Gal(K_S/K) \To Aut(\pi_1(\bar X)^{un}_{/\Q_l}) \To
Aut(H_1(\mathcal P))
$$
where $\mathcal P$ denote the nilpotent Lie algebra associated to
$\pi_1(\bar X)^{un}_{/\Q_l}$. This time, the Schlessinger
criteria may not help us in finding a universal representation.

In the first part of this paper, we have introduced a derivation
version [Tsu] which is a Schlessinger friendly representation of
$\Z_l$-Lie algebras

$$
\gal (K_S/K) \To Der(\mathcal P)/Inn(\mathcal P),
$$
or one could deform the following morphism, fixing its mod-$l$
reduction
$$
\gal (K_S/K) \To \out (\pi_1(\bar X)^{(l)}).
$$

\section{Deformation of representations of Lie algebras}

The action of $Aut(\pi_1(\bar X))$ on $H^i(\pi_1(\bar X),\Z_l)$
is compatible with the nondegenerate alternating form defined by
the cup product
$$
H^1(\pi_1(\bar X),\Z_l)\times H^1(\pi_1(\bar X),\Z_l)\To
H^i(\pi_1(\bar X),\Z_l)\cong \Z_l
$$
which is why the grade zero part of $Aut(\pi_1(\bar X))$ is the
same as $Sp(2g,\Z_l)$. As a Galois module, this is exactly the
Galois representation associated to the Tate module of the
Jacobian variety of $X$. This representation completely
determines the map
$$
\rho^{ni,l}_X:Gal(K_S/K) \To Aut(\mathcal U_S)
$$
in view of the isomorphism $\mathcal U_S \cong \out (\pi_1(\bar
X)^{(l)})\otimes \Q_l$. This is a convenient framework to perform
different versions of restricted deformation theories appearing
in the proof of Wiles using the language of Lie algebras. Our
ultimate goal is to make a Lie algebra version of computations in
[Wil],[Tay-Wil] and [Bre-Con-Dia-Tay]. One also has explicit
information about the Mal\v{c}ev pro-nilpotent Lie algebra
associated to $\pi_1(\bar X)^{un}_{/\Q_l}$ [Pri] and one can
deform the Galois action on this Lie algebra.

Let us concentrate on deforming the graded Lie algebra version of
the outer Galois representation
$$
\gal (K_S/K) \To \out (\pi_1(\bar X)^{(l)}).
$$
We could fix the mod-$l$ representation, or fix restriction of
this representation to decomposition Lie algebra $\mathcal D_p$
at prime $p$, which is induced by the same filtration as
$Gal(K_S/K)$ on the decomposition group. For each prime $p$ of
$K$ we get a map
$$\mathcal D_p \To \gal (K_S/K)\To \out (\pi_1(\bar X)^{(l)}).
$$
\begin{thm}
For a graded $\Z_l$-Lie algebra $L$, let $D(L)$ be the set of
representations of $\gal(K_S/K)$ to $L$ which reduce to
$$
\bar {\rho}: \gal (K_S/K) \To \out (\pi_1(\bar X)^{(l)})/l\out
(\pi_1(\bar X)^{(l)})
$$
after reduction modulo $l$. Assume that $\gal(K_S/K)$ is a free
$\Z_l$ Lie algebra. Then, there exists a universal deformation
graded $\Z_l$-Lie algebra $L_{univ}$ and a universal
representation
$$
\gal (K_S/K) \To L_{univ}
$$
representing the functor $D$. In case $\gal(K_S/K)$ is not free,
then one can find a hull for the functor $D$.
\end{thm}
\textbf {Proof.} For free $\gal (K_S/K)$ by theorem 2.10 in [Ras]
we have
$$
H^2(\gal (K_S/K), Ad\circ \bar{\rho} )=0
$$
which implies that $D$ is pro-representable. In case $\gal
(K_S/K)$ is not free, we have constructed a miniversal
deformation Lie algebra for another functor in theorem 2.11 of
[Ras], which implies that the first three Schlessinger criteria
hold. By a similar argument one could prove that there exists a
hull for $D$. Note that $\out (\pi_1(\bar X)^{(l)})$ is
pronilpotent, and for deformation of such an object one should
deform the truncated object and then take a limit to obtain a
universal object. $\Box$

\section{Grothendieck's anabelian conjectures}

The general philosophy of Grothendieck's anabelian conjectures is
to characterize invariants of a smooth hyperbolic curve by its
non-abelian geometric fundamental group. In particular,
Galois-module structure of the outer automorphism group of a
smooth hyperbolic curve defined over a number field, should
contain all the arithmetic information about $X$. This could not
be true if we consider the pro-$l$ completion of the fundamental
group. So, translating this to the language of Lie algebras, one
should work with adelic Lie algebras. On the other hand, if one
considers the action of Galois group by conjugation on $\out
(\pi_1(\bar X)^{(l)})$ for all $l$, This can be determined by the
set of $l$-adic representations induced by the Jacobian variety,
which gives us less than what we want.

We propose to consider the Galois action on different
realizations of the motivic fundamental group of a smooth
hyperbolic curve. This way we get NGLAs together with Galois
actions. It is natural to expect this to carry arithmetic
structure of $X$. Here is a version of "Hom" conjecture
\begin{conj}
Let $X$ and $X'$ denote smooth hyperbolic curves defined over
number field $K$. There exists a natural one-to-one correspondence
$$
Hom_{K}(X,X')\To Hom_{Gal(\bar {K}/K)}(\pi_1^{mot}(\bar X)
,\pi_1^{mot}(\bar X '))
$$
where $\pi_1^{mot}(\bar X)$ denotes the motivic fundamental group
of the curves over $\bar K$ and $Hom_{Gal(\bar {K}/K)}$ denotes
the set of Galois equivariant homomorphisms between the two NGLAs.
\end{conj}

Practically, this means that, by considering all realizations of
the motivic fundamental group of a hyperbolic curve, we can grasp
all arithmetic information encoded in the curve. We also
conjecture that there should be a motivic Galois group which make
the NGLA version of the above conjecture true.
\begin{conj}
There exists a motivic Galois pro-NGLA $\gal ^{mot}(\bar {K}/K)$
which canonically maps to $\pi_1^{mot}(\bar X)$ for all hyperbolic
smooth curve $X$ defined over the number field $K$ such that the
following map is a natural one-to-one correspondence
$$
Hom_{K}(X,X')\To Hom_{\gal ^{mot}(\bar {K}/K)}(\pi_1^{mot}(\bar X)
,\pi_1^{mot}(\bar X ')).
$$
\end{conj}

Note that, by a conjecture of Deligne, the $Z_l$-Lie algebras
$\gal (\bar {K}/K)$ associated to $l$-adic central series
filtrations on the pro-unipotent fundamental group of
thri-punctured sphere $\pi_1^{un}(\mathbb P^1\setminus
0,1,\infty)$ are induced from a single $\Z$-Lie algebra by
extension of schalars.

\section*{Acknowledgements}
I would like to thank O. Gabber, M. Kontsevich, R. Mikhailov, J.
Tilouine, S. Wewers for enjoyable conversations. Also, I wish to
thank institute des hautes \`{e}tudes scientifiques for warm
hospitality during which part of this work was written down.

Sharif University of Technology, e-mail: rastegar@sharif.ir

\end{document}